\documentclass [10pt,twoside,b5paper]{article}
\usepackage{amsfonts}
\usepackage{amsthm}
\usepackage{amsmath}
\usepackage{amstext}
\usepackage{amssymb}
\usepackage{mathrsfs}
\usepackage{epsf}              
\usepackage{graphicx}          
\usepackage{fancybox}          
\usepackage{color}             
\usepackage{fancyhdr}
\usepackage[hang,footnotesize]{caption2}  
\def
\mathcal{\mathscr}
\newfont{\aaa}{cmb10 at 18pt}
\newfont{\bbb}{cmb10 at 10pt}

\newcommand{\beq}{\begin{equation}}
\newcommand{\eeq}{\end{equation}}
\newcommand{\bey}{\begin{eqnarray}}
\newcommand{\eey}{\end{eqnarray}}
\newcommand{\beyy}{\begin{eqnarray*}}
\newcommand{\eeyy}{\end{eqnarray*}}

\begin{document}

\setcounter{page}{1}
\qquad\\[8mm]

\noindent{\aaa{Stochastic integration for fractional L\'{e}vy
process and stochastic differential equation driven by fractional
L\'{e}vy noise
$^{^{^{\displaystyle*}}}$}}\\[1mm]

\noindent{\bbb Xuebin L\"{u}$^{1,2}$, Wanyang Dai$^{1}$}\\[-1mm]

\noindent\footnotesize{1 Department of Mathematics, Nanjing
University , Nanjing, P. R. China 210093 \\
 2 College of Science,
Nanjing University of
Technology, Nanjing, P. R. China 210009 \\}\\[6mm]

\normalsize\noindent{\bbb Abstract}\quad  In this paper, based on
the white noise analysis of square integrable pure-jump L\'{e}vy
process given by \cite{Lokka}, we define the formal derivative of fractional L\'{e}vy
process  defined by the square
integrable pure-jump L\'{e}vy process as the fractional L\'{e}vy noises by considering fractional L\'{e}vy process as the generalized
functional of L\'{e}vy process, and then we define the Skorohod integral with
respect to the fractional L\'{e}vy process.
Moreover, we propose a class of stochastic
Volterra equations driven by fractional L\'{e}vy noises  and investigate the existence and uniqueness of their solutions; In addition, we   propose a class of stochastic
differential equations driven by fractional L\'{e}vy noises and prove that under the
Lipschtz and linear conditions there exists  unique stochastic distribution-valued  solution.
\vspace{0.3cm}

\noindent{\bbb Keywords}\quad White noise;  fractional L\'{e}vy
processes; stochastic ordinary
linear differential equations; stochastic Volterra equation \\
{\bbb MSC}\quad 60E07, 60G20, 60G51, 60G52, 60H40 \\[0.4cm]

\noindent{\bbb{1\quad Introduction}}\\[0.1cm]
\setcounter{equation}{0} \noindent  The study on fractional
processes started from the fractional
 Brownian motion     introduced by Kolmogrov
  \cite{Kol1} and popularized by Mandelbrot and Van Ness \cite{Mandelbrot}.  The
  self-similarity and long-range dependence properties make fractional Brownian
  motion
 suitable to model driving noise in different applications such as
 mathematical finance and network traffic analysis. However, its
 light tails are often inadequate to model the higher variability
 phenomena appeared in these practical systems.

\vspace{-1mm}
\noindent \hrulefill\hspace{117mm}\\
{\footnotesize Corresponding author: Xuebin L\"{u},
E-mail: lvxuebin2008@163.com  \\
Project   Supported by Natural Science Foundation of China with
granted No. 10971249 and 11001051,10971076.}
\newpage
\noindent  Thus, it is  natural to consider the more general
fractional
 processes. Marquardt  \cite{Marquardt} introduced the fractional L\'{e}vy processes,  restricted to  the case of
 L\'{e}vy processes  with zero mean, finite variance and without
Brownian components. In \cite{HuangLC} and \cite{HuangLP} by white
noise approach, the authors constructed generalized fractional
L\'{e}vy processes as L\'{e}vy white noise functionals under a
simple condition on L\'{e}vy measure. In our previous work
\cite{HuangLU} \cite{LUHuang}, we defined infinite-dimensional
fractional L\'{e}vy processes on Gel$^{\prime}$fand Triple and
investigate its properties of distribution and sample path.

\par In order to use the fractional L\'{e}vy processes to
model the higher variability phenomena in real-world systems, it is
necessary to investigate stochastic calculus  for  fractional
L\'{e}vy processes and stochastic differential equations driven by
these processes. In \cite{Marquardt}, the authors defined stochastic
integral for a class of deterministic integrands with respect to
real-valued  fractional L\'{e}vy processes; In \cite{LUHuang}, we
defined stochastic integral for a class of real deterministic
functions and deterministic operator-valued processes with respect
to fractional L\'{e}vy processes on Gel$^{\prime}$fand Triple; In
\cite{Bender}, by using S-transform, the authors constructed the
stochastic calculus for convoluted L\'{e}vy processes which are
built by convoluting a Volterra-type kernel with a pure jump, zero
expectation L\'{e}vy process with finite moment of any order.
Especially, the authors investigate the Skorohod integral for
fractional L\'{e}vy processes whose underlying L\'{e}vy process has
finite moment of any order.
\par
The purpose of this paper is to investigate stochastic calculus of
fractional L\'{e}vy processes whose underlying  L\'{e}vy processes
are square integrable and investigate stochastic differential
equation driven by fractional L\'{e}vy processes.
\par
Motivated by the study of Lokka and Proske \cite{Lokka} about  the
infinite dimensional analysis of the square integrable pure jump
L\'{e}vy process on the Poisson space,  we investigate the
stochastic calculus for fractional L\'{e}vy processes whose
underlying L\'{e}vy proceses are only squre integrable. Our results
generalize that of \cite{Bender} which demands the underlying
L\'{e}vy processes have finite moment of any order.  Moreover, we
investigate stochastic Volterra equation driven by fractional
L\'{e}vy noise. In addition, we obtain a unique local continuous
solution for the stochastic differential equation driven by
fractional L\'{e}vy noise with Lipschtz and linear conditions.
\par
 This paper is organized as follows: In Section 2, we
recall the basic results about the infinite dimensional analysis of
the square integrable pure jump L\'{e}vy process given by  \cite{Lokka}, all of our work is done in this
framework; Based on the infinite dimensional analysis of L\'{e}vy
process, in Section 3, we define Skorohod integral with respect to
fractional L\'{e}vy processes and investigate its S-transform.
Moreover, we obtain an integration transformation formula which can
be used to
 transform stochastic differential equations driven by
 fractional L\'{e}vy  noises with different parameter into
 the same kind of stochastic differential equations with the same
 parameter  $\beta$;   In Section 4, we investigate the condition of existence
and uniqueness of  the solution  for stochastic differential equation
of Volterra type driven by  fractional L\'{e}vy processes. In
Section 5, we obtain a unique continuous global  solution for a
stochastic differential equation driven by fractional L\'{e}vy noise
with Lipschtz and linear conditions.
\\[4mm]
\noindent{\bbb 2\quad Infinite dimensional analysis of
square integrable pure jump L\'{e}vy process}\\[0.1cm]
 The    L\'{e}vy process $X=\{X_{t},t\geq0\}$
defined on a probability space $(\Omega,\mathcal{F}$,$P)$   is a
stochastic continuous process with stationary and independent
increments and the characteristic function of $X_{1}$ takes the
form:
$$
\widehat{\mu}(u)=\exp\{i\gamma
u-\frac{1}{2}u^{2}\sigma^{2}+\int_{\mathbb{R}}[e^{iux}-1-iux1_{|x|\leq
1}]d\nu(x)\},u \in \mathbb{R}, \eqno(2.1)
$$
where $\gamma\in\mathbb{R}$, $\sigma\geq 0$, and $\nu$ is the
L\'{e}vy
 measure  satisfying  $\nu(\{0\})=0$ and
$$
 \int_{\mathbb{R}}(|x|^{2}\wedge 1 )d\nu(x)<\infty.\eqno(2.2)
$$
 Then, for any $t\geq0$, $X$ has the L\'{e}vy-It\^{o} decomposition
$$
 X_{t}=\gamma t+B_{t}+\int_{0}^{t}\int_{|x|\geq1}xN(ds,dx)+
\int_{0}^{t}\int_{|x|\leq1}x\widetilde{N}(ds,dx), \eqno(2.3)
 $$
 where
 $$
 N((0,t]\times A)=\sum_{s\leq t}1_{A}(\triangle
X_{s}), A\in \mathcal{B}
 (\mathbb{R})
 $$
is a Poisson random measure on $\mathbb{R}_{+}\times \mathbb{R}_{0}$
$(\mathbb{R}_{0}=:\mathbb{R}\setminus\{0\})$  and
$$
 \widetilde{N}((0,t]\times A)=N((0,t]\times A)-t\nu(A) \eqno(2.4)
$$
is its compensation. $B=\{B_{t},t\geq0\}$ is a Brownian motion with
mean 0 and covariance operator $\sigma^{2}$ which is independent of
$N$.
 The first integral in (2.3) converges in probability (even a.s.),
the second one converges in $ L^{2}(\Omega, \mathcal{F},P)$. If in
(2.1),   $\sigma=0$, we call $X$ a  L\'{e}vy process without
Brownian component. In this paper, we always assume that the
L\'{e}vy process $X$ has no Brownian part. Furthermore, we suppose
$\mathbb{E}[X_{1}]=0$  and
$$
 \int_{|x|>1}|x|^{2}d\nu(x)<\infty.
$$
Thus, (2.1) can be written as
$$
\widehat{\mu}(u)=\exp\{\int_{\mathbb{R}}[e^{iux}-1-iux]d\nu(x)\},u
\in \mathbb{R},
$$
and
$$
X_{t}=\int_{0}^{t}\int_{\mathbb{R}}x\widetilde{N}(ds,dx).
$$
In this case, $X$ is a martingale and we call it pure-jump L\'{e}vy
process. We will work with a two-side L\'{e}vy
process constructed by taking two independent  copies $X^{(1)}=\{X_{t}^{(1)},t\geq0\}$
and $X^{(2)}=\{X_{t}^{(2)},t\geq0\}$ of a one-side L\'{e}vy
process and setting
\begin{equation}
X_{t}=\begin{cases}X_{t}^{(1)},&\  t\geq0\\ -X_{-t_{-}}^{(2)},&\ t<0
\ .\nonumber
\end{cases}
\end{equation}
\par Next we recall the basic results about the infinite
dimensional analysis of the square integrable pure jump L\'{e}vy
process given by \cite{Lokka}.
\par
 Let $\xi_{n}$ denote denote the
$n'$th Hermite function, the set of Hermite functions
$\{\xi_{n}\}_{n\in\mathbb{N}}$ is an orthonormal basis of
$L^{2}(\mathbb{R})$. Denote by $\mathcal{S}$$(\mathbb{R}^{d})$ the
Schwartz space of rapidly decreasing $C^{\infty}$-functions on
$\mathbb{R}^{d}$ and by $\mathcal{S}'$$(\mathbb{R}^{d})$ the space
of tempered distributions. The nuclear topology on
$\mathcal{S}$$(\mathbb{R}^{d})$ is induced by the pre-Hilbertian
norms
$$
\|\phi\|_{p}^{2}:=\sum_{\alpha=(\alpha_{1},
\ldots,\alpha_{d})\in\mathbb{N}^{d}
  }(1+\alpha)^{2p}(\phi,\xi_{\alpha})_{L^{2}(\mathbb{R}^{d})}^{2},p\in\mathbb{N}
  _{0},
$$
where $(1+\alpha)^{2p}=\prod_{i=1}^{d}(1+\alpha_{i})^{2p}$,
$\xi_{\alpha}(x_{1},
\ldots,x_{d})=\prod_{i=1}^{d}\xi_{\alpha_{i}}(x_{i}) $.
 Let $\mathbb{U}=\mathbb{R}\times \mathbb{R}_{0}$,   define
$$
\mathcal{S}(\mathbb{U}):=\{\phi\in
\mathcal{S}(\mathbb{R}^{2}):\phi(x,0)=\frac{\partial\phi}{\partial
y}(x,0)=0 \}.
$$
$\mathcal{S}(\mathbb{U})$ is a closed subspace of
$\mathcal{S}(\mathbb{R}^{2})$, thus it is a countably Hilbertian
nuclear algebra endowed with the topology induced by the norms  $\|
\cdot\|_{p}$, and its dual
$\mathcal{S}'(\mathbb{U})\supset\mathcal{S}'(\mathbb{R}^{2})$. For
$\phi\in\mathcal{S}(\mathbb{U})$, $\Phi\in\mathcal{S}'(\mathbb{U})$,
the action of $\Phi$ on $\phi$ is given by
 $
 \langle\Phi,\phi\rangle=\int_{\mathbb{U}}\Phi(x,y)\phi(x,y)dxdy.
 $
Assume that $\nu$ is the L\'{e}vy
 measure  on $\mathbb{R}_{0}$ satisfying
$$
 \int_{\mathbb{R}_{0}}|x|^{2}d\nu(x)<\infty.\eqno(2.5)
$$
Denote $\lambda$ the Lebesgue measure on $\mathbb{R}$ and let $\pi$
denote the measure on $\mathbb{U}$ given by $\pi=\lambda\times\nu$.
By Lemma 2.1 of \cite{Lokka}, there exists an element denoted by
$1\otimes\dot{\nu}$ in $\mathcal{S}'(\mathbb{U})$ such that
$$
\langle 1\otimes\dot{\nu},
\phi\rangle=\int_{\mathbb{U}}\phi(x)\pi(dx),\phi\in\mathcal{S}(\mathbb{U}).\eqno(2.6)
$$
Denote ${L}^{2}(\mathbb{U} ,\pi )$ by the space of all square
integrable functions on $\mathbb{U} $ with respect to $\pi $, let
$(\cdot,\cdot)_{\pi}$ the inner product on ${L}^{2}(\mathbb{U},\pi)$
and $|\cdot|_{\pi}$ the corresponding norms on this space. Define  $
\mathcal{N}_{\pi}:=\{\phi\in\mathcal{S}(\mathbb{U}):|\phi|_{\pi}=0\}
 $,  then $\mathcal{N}_{\pi}$ is a closed ideal of
$\mathcal{S}(\mathbb{U})$.
 Let $\widetilde{\mathcal{S}}(\mathbb{U})$ be
the space $
\widetilde{\mathcal{S}}(\mathbb{U}):=\mathcal{S}(\mathbb{U})/\mathcal{N}_{\pi}
$ endowed with the topology induced by the system of norms $
\|\widehat{\phi}\|_{p,\pi}:=\inf_{\psi\in\mathcal{N}_{\pi}}\|\phi+\psi\|_{p},
$ then $\widetilde{\mathcal{S}}(\mathbb{U})$ is a nuclear algebra.
Let $\widetilde{\mathcal{S}}'(\mathbb{U})$ be the dual of
$\widetilde{\mathcal{S}}(\mathbb{U})$, and for $p\in\mathbb{N}$, let
$\widetilde{\mathcal{S}}_{p}(\mathbb{U})$ denote the completion of
$\widetilde{\mathcal{S}}(\mathbb{U})$ with respect to the norm
$\|\cdot\|_{p,\pi}$, $\widetilde{\mathcal{S}}'_{-p}(\mathbb{U})$
denote the dual of $\widetilde{\mathcal{S}}_{p}(\mathbb{U})$.
$\widetilde{\mathcal{S}}(\mathbb{U})$ is the projective limit of
$\{\widetilde{\mathcal{S}}_{p}(\mathbb{U}), p\geq 0\}$, and
$\widetilde{\mathcal{S}}'(\mathbb{U})$ is the inductive limit of
$\{\widetilde{\mathcal{S}}'_{-p}(\mathbb{U}), p\geq 0\}$.
$\widetilde{\mathcal{S}}(\mathbb{U})$ has similar nice properties as
the classical Schwartz space. Thus, Lokka and Proske \cite{Lokka}
introduce it to   construct   the white noise analysis of L\'{e}vy
process. \\
  \textbf{Theorem 2.1} (\cite{Lokka})\quad(1)There exists a probability measure $\mu_{\pi}$ on
$\widetilde{\mathcal{S}}'(\mathbb{U})$ such that
$$
\int_{\widetilde{\mathcal{S}}'(\mathbb{U})}e^{i\langle \omega,
\phi\rangle}d\mu_{\pi}(\omega)=\exp\{\int_{\mathbb{U}}(e^{i\phi(x)}-1)d\pi(x)\},
\forall\phi\in\widetilde{\mathcal{S}}(\mathbb{U}).\eqno(2.7)
$$
(2) There exists a $p_{0}\in\mathbb{N}$ such that
$1\otimes\dot{\nu}\in
\widetilde{\mathcal{S}}_{-p_{0}}'(\mathbb{U})$, and a natural number
$q_{0}>p_{0}$ such that the imbedding operator
$\widetilde{\mathcal{S}}_{-p_{0}}'(\mathbb{U})\hookrightarrow\widetilde{\mathcal{S}}_{-q_{0}}'(\mathbb{U})$
is Hilbert-Schimidt and
$\mu_{\pi}(\widetilde{\mathcal{S}}_{-q_{0}}'(\mathbb{U}))=1$.
\par
From now on, for all $q_{0}$, $p_{0}$ are described in the Theorem
2.1. Set $\Omega=\widetilde{\mathcal{S}}'(\mathbb{U})$ and
$\mu_{\pi}$ given by Theorem 2.1 on which \cite{Lokka} give the
infinite dimensional calculus for pure jump measure, and all of our
following discussion is based on this probability space.
\par
Let $C_{n}(\cdot)$ be the generalized Charlier polynomials given by
\cite{Lokka},  for all $m, n\in \mathbb{N}$,
$\varphi^{(n)}\in\widetilde{\mathcal{S}}(\mathbb{U})^{\widehat{\otimes}n}$,
$\psi^{(m)}\in
\widetilde{\mathcal{S}}(\mathbb{U})^{\widehat{\otimes}m}$,
($\widehat{\otimes}$  denotes the symmetrized tensor product),  the
following orthogonality relation holds,
\begin{equation}
\int_{\widetilde{\mathcal{S}}'(\mathbb{U})}\langle C_{n}(\omega),
\varphi^{(n)}\rangle\langle
C_{m}(\omega), \psi^{(m)}\rangle d\mu_{\pi}(\omega)=\begin{cases}0,&\  n\neq m\\
n!(\varphi^{(n)},\psi^{(n)})_{\pi},&\ n=m \ .\nonumber
\end{cases}
\end{equation}
Especially, for $n=1$, $C_{1}(\omega)=\omega- 1\otimes\dot{\nu}$.
Since  $L ^{2}(\mathbb{U})$ is dense in
$\widetilde{\mathcal{S}}(\mathbb{U})$, for $f  \in L
^{2}(\mathbb{U})$, there exists   a sequence of functions $f _{n}\in
\widetilde{\mathcal{S}}(\mathbb{U}) $ such that $f_{n}\rightarrow f
$ in $\widehat{L}^{2}(\mathbb{U} ,\pi) $ as $n\rightarrow\infty$.
Define $\langle C_{1}(\omega), f\rangle$ by $ \langle C_{1}(\omega),
f \rangle=\lim_{n\rightarrow\infty}\langle C_{1}(\omega),
f_{n}\rangle$ $  (limit \ in \ L^{2}(\mu_{\pi}))$, the definition is
independent of the choice of approximating sequence  and   the
following isometry holds
$$
\int_{\widetilde{\mathcal{S}}'(\mathbb{U})}\langle
C_{1}(\omega),f\rangle^{2}
d\mu_{\pi}(\omega)=\int_{\widetilde{\mathcal{S}}'(\mathbb{U})}\langle
\omega-1\otimes\dot{\nu},f\rangle^{2}
d\mu_{\pi}(\omega)=|f|_{\pi}^{2}.\eqno(2.8)
$$
For any Borel sets $\Lambda_{1}\subset \mathbb{R}$ and
$\Lambda_{2}\subset \mathbb{R}_{0}$ such that the 0 is not in the
closure of $\Lambda_{2}$, define the random measure
$$
N(\Lambda_{1},\Lambda_{2}):=\langle
\omega,1_{\Lambda_{1}\times\Lambda_{2}}\rangle,
\widetilde{N}(\Lambda_{1},\Lambda_{2}):=\langle
\omega-1\otimes\dot{\nu},1_{\Lambda_{1}\times\Lambda_{2}}\rangle.
$$
From the characterization function of $\mu_{\pi}$, it is easy to
deduced that  $N$ is a Poisson random measure, and $\widetilde{N}$
is the corresponding compensated measure. The compensator of $
N(\Lambda_{1},\Lambda_{2})$ is given by $\langle
1\otimes\dot{\nu},1_{\Lambda_{1}\times\Lambda_{2}}\rangle $ which is
equal to $\pi(\Lambda_{1}\times\Lambda_{2})$. Moreover,
$$
\int_{\mathbb{U}}\phi(s,x)\widetilde{N}(ds,dx)=\langle
\omega-1\otimes\dot{\nu},\phi\rangle, \phi\in
L^{2}(\mathbb{U},\pi).\eqno(2.9)
$$
Then, the pure jump L\'{e}vy process has a representation
$$
X_{t}=\int_{0}^{t}\int_{\mathbb{R}}x\widetilde{N}(ds,dx)=\langle
\omega-1\otimes\dot{\nu},f_{t}\rangle,f_{t}(s,x)=x1_{s\leq
t}(s).\eqno(2.10)
$$
The Wiener integral of $g\in L^{2}(\mathbb{R})$ with respect to $X$
has the following representation
 $$
\int_{\mathbb{R}}g(s)dX_{s}=\int_{\mathbb{R}}\int_{\mathbb{R}}g(s)x\widetilde{N}(ds,dx)
=\langle \omega-1\otimes\dot{\nu},\phi\rangle,\eqno(2.11)
$$
where $ \phi(s,x)=xg(s)$.
\par
 Define the space $
\mathcal{P}(\widetilde{\mathcal{S}}'(\mathbb{U})) =
\{f:\widetilde{\mathcal{S}}'(\mathbb{U})\rightarrow C,
f(\omega)=\sum_{n=0}^{N}\langle \omega^{\otimes n},
\phi^{(n)}\rangle,$     $\omega\in\widetilde{
\mathcal{S}}'(\mathbb{U}), \phi^{(n)}\in\widetilde{
\mathcal{S}}(\mathbb{U})^{\widehat{\otimes}n}, N\in\mathbb{N}\}
 $,
 $ f$ is called a continuous polynomial function if
$f\in\mathcal{P}(\widetilde{\mathcal{S}}'(\mathbb{U}))$  and it
admits a unique representation of the form
$$
f(\omega)=\sum_{n=0}^{\infty}\langle C_{n}(\omega),
f_{n}\rangle,f_{n}\in\widetilde{
\mathcal{S}}(\mathbb{U})^{\widehat{\otimes}n}.
$$
For any number  $p\geq q_{0}$, define the Hilbert space
$(\mathcal{S})_{p}^{1}$ as the completion of
$\mathcal{P}(\widetilde{\mathcal{S}}'(\mathbb{U}))$ with respect to
the norm
$$
\|f\|^{2}_{p,1}=\sum_{n=0}^{\infty}(n!)^{2}\|f_{n}\|_{p,\pi}^{2}.
$$
The corresponding inner product is
$$
((f,g))_{p,1}=\sum_{n=0}^{\infty}(n!)^{2}((f_{n},g_{n}))_{p,\pi}.
 $$
where $((\cdot,\cdot))_{p,\pi}$ denotes the inner product on
$\widetilde{\mathcal{S}}_{p}(\mathbb{U})^{\widehat{\otimes}n}$.
Obviously, $(\mathcal{S})_{p+1}^{1}\subset(\mathcal{S})_{p}^{1}$. In
\cite{Lokka}, the authors define $(\mathcal{S})^{1}$ as the
projective limit of $\{(\mathcal{S})_{p}^{1}, p\geq q_{0}\}$, and it
is a nuclear Fr\'{e}chet space which can be densely imbedding in
$L^{2}(\mu_{\pi})$.  Denote $(\mathcal{S})_{-p}^{-1}$ as the dual of
$(\mathcal{S})_{p}^{1}$, $(\mathcal{S})^{-1}$ as the inductive limit
of $\{(\mathcal{S})_{-p}^{-1}, p\geq q_{0}\}$ which is equal to the
dual of  $(\mathcal{S})^{1}$. $F\in(\mathcal{S})^{-1}$ if and only
if $F$ admits an expansion
$$
F(\omega)=\sum_{n=0}^{\infty}\langle C_{n}(\omega),
F_{n}\rangle,F_{n}\in\widetilde{
\mathcal{S}}'(\mathbb{U})^{\widehat{\otimes}n},
$$
and there exists a $ p\geq q_{0}$ such that
$$
\|F\|^{2}_{-p,-1}=\sum_{n=0}^{\infty}\|F_{n}\|_{-p,\pi}^{2}<\infty .
$$
For  $F\in(\mathcal{S})^{-1} $,  $f\in(\mathcal{S})^{1}$, $$
\langle\langle F,g\rangle\rangle=\sum_{n=0}^{\infty}n!\langle
F_{n},f_{n}\rangle_{\pi},
$$
where $ \langle\langle \cdot,\cdot\rangle\rangle$ is an extension of
the inner product on $L^{2}(\mu_{\pi})$. $(\mathcal{S})^{1}$ is
called space of stochastic test functions, $(\mathcal{S})^{-1}$ is
called space of stochastic distribution functions, they are pairs of
dual spaces, $(\mathcal{S})^{1}\subset
L^{2}(\mu_{\pi})\subset(\mathcal{S})^{-1}$.
\par
Next we recall the S-transform given by \cite{Lokka} which can
transform stochastic distribution functions to deterministic
functionals. Let
$$
\widetilde{e}(\phi,\omega):=\exp(\langle\omega,\ln(1+\phi)\rangle-\langle1\otimes\dot{\nu},\phi\rangle),
$$
it is analytic as a function of
$\phi\in\widetilde{\mathcal{S}}_{q_{0}}$ for
$\phi\in\widetilde{\mathcal{S}}_{q_{0}}$ satisfying $\phi(x)>-1$ for
all $x\in \mathbb{U}$. Moreover, it has the following chaos
expansion,
$$
\widetilde{e}(\phi,\omega)=\sum_{n=0}^{\infty}\frac{1}{n!}\langle
C_{n}(\omega), \phi^{\otimes n}\rangle.
$$
  Denote
 $ U_{p}:=\{\phi\in\widetilde{\mathcal{S}}(\mathbb{U}):
\|\phi\|_{p,\pi}<1\}  $, by the chaos expansion of
$\widetilde{e}(\phi,\omega)$
$\widetilde{e}(\phi,\omega)\in(\mathcal{S})^{1}_{p}$ if and only if
$\phi\in U_{p}$.\\
 \textbf{Definition 2.2}(\cite{Lokka})\ Let $F\in(\mathcal{S})^{-1}_{-p} $, $\xi\in
 U_{p}$, the the S-transform of $F$ is defined by
 $$
 S(F)(\xi):=\langle\langle
 F,\widetilde{e}(\xi,\omega)\rangle\rangle.
$$
 For example, if $F=\sum_{n=0}^{\infty}\langle
C_{n}(\omega),F_{n}\rangle\in(\mathcal{S})^{-1}_{-p} $, $\xi\in
 U_{p},$ then  $
 S(F)(\xi)=\sum_{n=0}^{\infty}\langle F_{n},  \xi^{\otimes
 n}\rangle_{\pi}.
 $ \par
 Denote $\mathcal{U}=Hol(0)$  the algebra of germs of functions
that are holomorphic in a neighborhood of 0. The S-transform is
isomorphic between $(\mathcal{S})^{-1} $ and $\mathcal{U}$.
  Since  $f, g\in\mathcal{U}$, then $fg\in\mathcal{U}$, then   the following definition of  Wick product is well-defined.\\
 \textbf{Definition 2.3}(\cite{Lokka})
(Wick product) \ Let $F,G\in (\mathcal{S})^{-1}$, define the Wick
product $F\diamond G$ of $F$ and $G$ by
$$
F\diamond G=S^{-1}( S(F)S(G)).
$$
  Let $F:\mathbb{U}\rightarrow (\mathcal{S})^{-1}$ be
the random fields with chaos expansion
$$
F(x)=\sum_{n=0}^{\infty}\langle C_{n}(\omega),F_{n}(\cdot,x)\rangle,
$$
where $
F_{n}(\cdot,x)\in\mathcal{S}'(\mathbb{U})^{\widehat{\otimes}n}$ and
$\|F(x)\|_{-p,-1}<\infty, $ for some $p>0$. Let $\mathbb{L}$ denote
the set of all $F:\mathbb{U}\rightarrow (\mathcal{S})^{-1}$ such
that $
\widetilde{F_{n}}\in\mathcal{S}'(\mathbb{U})^{\widehat{\otimes}(n+1)}$($
\widetilde{F_{n}}$ is the symmetrization of $F$) and $
\sum_{n=0}^{\infty}|\widetilde{F_{n}}|_{-p,\pi}^{2}<\infty  $
for some $p>0$.\\
 \textbf{Definition 2.4}(\cite{Lokka})(Skorohod integral ) For
 $F\in\mathbb{L}$, define the Skorohod integral $\delta(F)$ by
 $$
 \delta(F):=\sum_{n=0}^{\infty}\langle
 C_{n+1}(\omega),\widetilde{F_{n}}\rangle.
$$
From the assumption on $\mathbb{L}$,  we see that
$\delta(F)\in(\mathcal{S})^{-1}$. For the predictable integrands,
the Skorohod integral coincides with the usual Ito-type integral
with respect to the compensated Poisson random measure.  \\
 \textbf{Proposition 2.5}(\cite{Lokka}) If $F\in\mathbb{L}$, then
$\delta(F)\in(\mathcal{S})^{-1}_{-p}$ for some $p>0$ and
$$
S\delta(F)(\xi)=\int_{\mathbb{U}} S F(x)(\xi)\xi(x)\pi(dx), \xi\in
U_{p}.\eqno(2.12)
$$
\\[4mm]
\noindent{\bbb 3\quad  S-transform and  Skorohod integral for fractional  L\'{e}vy processes }\\[0.1cm]
\noindent In this section, we give the S-transform and the Skorohod
integral for fractional  L\'{e}vy processes  based on Section 2.1.
\par
First, we recall the definition of fractional L\'{e}vy processes
(for more details, see \cite{Marquardt}, \cite{HuangLU},
\cite{LUHuang}). The $\beta$-fractional L\'{e}vy process
$\{X_{t}^{\beta},t\geq0\}$ ($0<\beta<\frac{1}{2})$  is defined by:
$$
X_{t}^{\beta}=\int_{\mathbb{R}}I_{-}^{\beta}\chi_{[0,t]}(s)dX_{s}=\frac{1}{\Gamma(\beta+1)}
\int_{-\infty}^{\infty}((t-s)_{+}^{\beta}-(-s)_{+}^{\beta})dX_{s},\
\eqno(3.1)
$$
where $X$ is a two-side L\'{e}vy process satisfying all the
assumptions in Section 2,\begin{equation}
\chi_{[0,t]}(s)=\begin{cases}1,& 0<s<t\\ -1,&\ t<s<0 \\0, &\
\text{else}.\nonumber
\end{cases}
\end{equation}
Furthermore, $I_{-}^{\beta}$ is the Riemann-Liouville fractional
integral operator defined by
$$
(I_{-}^{\beta}f)(t)=\frac{1}{\Gamma(\beta)}\int_{t}^{+\infty}(s-t)^{\beta-1}f(s)ds\
, f\in\mathcal{S}(\mathbb{R}),
$$
where $x_{+}=x\vee0$, and $\Gamma(\cdot)$ is the Gamma function.
\par
Define $(K^{\beta}f)(s,x):=xI_{-}^{\beta}f(s)$, the
Riemann-Liouville fractional differential operator are applied only
to the time variable $s$. Since $I_{-}^{\beta}\chi_{[0,t]}\in L^{2}$,
by (2.11), $ X_{t}^{\beta}$ has the following representation
$$
X_{t}^{\beta}=\langle C_{1},
K^{\beta}\chi_{[0,t]}\rangle=\delta(K^{\beta}\chi_{[0,t]}).\eqno(3.2) $$
Thus, by (2.11), we  get the S-transform of fractional L\'{e}vy process \\
$$
SX_{t}^{\beta}(\eta)=\int_{\mathbb{R}}\int_{\mathbb{R}_{0}}I_{-}^{\beta}\chi_{[0,t]}(s)y\eta(s,y
)\nu(dy)ds, \eta\in U_{p}, p> q_{0}. \eqno(3.3)  $$
 For $g$ satisfying $I_{-}^{\beta}g\in L^{2}$, the Wiener
integral with respect to $ X^{\beta}$ can be written as
$$
\int_{\mathbb{R}}g(t)dX_{t}^{\beta}=\int_{\mathbb{R}}I_{-}^{\beta}g(t)dX_{t}
 =\langle C_{1},
K^{\beta}g\rangle=\delta(K^{\beta}g).\eqno(3.4)
$$
Since $ \int_{\mathbb{R}}g(t)dX_{t}^{\beta}\in L^{2}(\Omega)$, its
S-transform is given by
$$
S(\int_{\mathbb{R}}g(t)dX_{t}^{\beta})(\eta)
=\int_{\mathbb{R}}\int_{\mathbb{R}_{0}}I_{-}^{\beta}g(t)y\eta(t,y)\nu(dy)dt,
\eta\in U_{p},  p> q_{0}.  \eqno(3.5)
$$
By the following fractional integral by parts formula of operator
$I_{\pm}^{\beta}$:
$$
\int_{-\infty}^{\infty}f(s)I_{+}^{\beta}g(s)ds=
\int_{-\infty}^{\infty}g(s)I_{-}^{\beta}f(s)ds, \ f,g\in
\mathcal{S}(\mathbb{R})
$$
which can be extended to $f\in L^{p}(\mathbb{R})$, $g\in
L^{r}(\mathbb{R})$ with $p>1$, $r>1$ and
$\frac{1}{p}+\frac{1}{r}=1+\beta$, where
$$
(I_{+}^{\beta}f)(t)=\frac{1}{\Gamma(\beta)}\int_{-\infty}^{t}(t-s)^{\beta-1}f(s)ds\
,
$$ (3.3) can be written as
$$
SX_{t}^{\beta}(\eta)=\int_{\mathbb{R}}\int_{\mathbb{R}_{0}}1_{[0,t]}(s)yI_{+}^{\beta}\eta(\cdot,y
)(s)\nu(dy)ds
=\int_{0}^{t}\int_{\mathbb{R}_{0}}yI_{+}^{\beta}\eta(\cdot,y
)(s)\nu(dy)ds,t\geq 0.
$$
Hence,
$$
\frac{d}{dt}SX_{t}^{\beta}(\eta)=\int_{\mathbb{R}_{0}}I_{+}^{\beta}\eta(\cdot,y
)(t)y\nu(dy), \eta\in U_{p}, p> q_{0}, t\geq 0. \eqno(3.6)
$$
Note that (3.6) also holds for $t<0$. We denote
$\dot{X}^{\beta}_{t}$  the fractional L\'{e}vy noise in the
following sense:
$$
S\dot{X}^{\beta}_{t}(\eta)=\frac{d}{dt}SX^{\beta}_{t}(\eta), \eta\in
U_{p},  p> q_{0} .
$$
Then we can prove that $\dot{X}^{\beta}_{t}$ is a generalized
stochastic distribution function and it has a chaos expansion.\\
\textbf{Theorem 3.1}\   $\dot{X}^{\beta}_{t}\in
(\mathcal{S})^{-1}_{-p}$ for all $p>\max\{1,q_{0}\}$  and
$$
\dot{X}^{\beta}_{t}=\langle C_{1}, \lambda_{t}\rangle, \eqno(3.7)
$$
where
$$\lambda_{t}(u,y)=yI_{-}^{\beta}\delta_{t}(u)=y(t-u)_{+}^{\beta-1}/\Gamma(\beta)$$
 \textbf{Proof:}  We first show that  $
\langle C_{1}, \lambda_{t}\rangle\in (\mathcal{S})^{-1}_{-p}$ for
all $p>\max\{1,q_{0}\}$. By the estimate
$$
|\int_{\mathbb{R}}(t-u)_{+}^{\beta-1} \xi_{n}(u)du |  \leq C
n^{\frac{2}{ 3}-\frac{\beta}{2}},
$$
from section 4 of \cite{Elliott} , where $C$ is a certain constant
independent of $t$,
$$
\aligned
\|\langle C_{1}, \lambda_{t}\rangle\|_{-1,-p}^{2}
&=\frac{\int_{\mathbb{R}}|y|^{2}d\nu(y)}{\Gamma(\beta)}\sum_{n=1}^{\infty}(n+1)^{-2p}\langle(t-\cdot)_{+}^{\beta-1},
\xi_{n}\rangle_{L^{2}(\mathbb{R})}^{2}\\
&=A\sum_{n=1}^{\infty}(n+1)^{-2p}
(\int_{\mathbb{R}}(t-u)_{+}^{\beta-1} \xi_{n}(u)du)^{2}\\
 &\leq AC\sum_{n=1}^{\infty}(n+1)^{-2p+\frac{4}{3}-\beta}\\
&< +\infty,  for\  p>\max\{1,q_{0}\},
\endaligned  \eqno(3.8)
$$
where
$A=\frac{\int_{\mathbb{R}}|y|^{2}d\nu(y)}{\Gamma(\overline{\beta})}$
is a positive constant. Thus, $ \langle C_{1}, \lambda_{t}\rangle\in
(\mathcal{S})^{-1}_{-p}$ for all $p>\max\{1,q_{0}\}$. Next, we prove
(3.7) holds. In fact,
$$
(I_{-}^{\beta}\delta_{t})(s)=\frac{1}{\Gamma(\beta) }
\int_{\mathbb{R}_{+}}\frac{\delta_{t}(s+u)du}{u
^{1-\beta}}=\frac{(t-s)_{+}^{\beta-1}}{\Gamma(\beta)}.
$$
Taking S-transform of $ \langle C_{1}, \lambda_{t}\rangle$, we have
$$\aligned
 S\langle C_{1}, \lambda_{t}\rangle(\eta)
&=\int_{\mathbb{R}}\int_{\mathbb{R}_{0}}\frac{y(t-s)_{+}^{\beta-1}}{\Gamma(\beta)}
\eta(s,y)\nu(dy)ds\\
&=\int_{\mathbb{R}}\int_{\mathbb{R}_{0}}yI_{-}^{\beta}\delta_{t}(s)\eta(s,y)\nu(dy)ds\\
&=\int_{\mathbb{R}}\int_{\mathbb{R}_{0}}y\delta_{t}(s)I_{+}^{\beta}\eta(\cdot,y)(s)\nu(dy)ds\\
&=\int_{\mathbb{R}_{0}}yI_{+}^{\beta}\eta(\cdot,y)(t)\nu(dy),
\eta\in U_{p},  p> q_{0}.\endaligned
$$
Hence, it follows from (3.6) that (3.7) holds.\hfill$\Box$

By (3.5) for $g$ satisfying $I_{-}^{\beta}g\in L^{2}$,  we have
  $$ \int_{\mathbb{R}}g(t)dX_{t}^{\beta}=\int_{\mathbb{R}}g(t)\diamond
\dot{X}^{\beta}_{t}dt.\eqno(3.9)
$$
Based on  (3.9) we can define Skorohod integral for
$(\mathcal{S})^{-1}$-valued processes with respect to the fractional
L\'{e}vy process $ X^{\beta}$ as follows. \\
 \textbf{Definition 3.2}(1)  $F: \mathbb{R}\rightarrow(\mathcal{S})^{-1}$ is differentiable,
 if $\forall t  \in \mathbb{R}$,
$$
 \lim_{\triangle t\rightarrow 0}\frac{ F(t+\triangle t) -
F(t) }{\triangle t}
  $$
 exists in  $(\mathcal{S})^{-1}$ .\\
 (2) Suppose  $F:\mathbb{R}\longrightarrow
 (\mathcal{S})^{-1}$ is a given function such that
 $
 \langle\langle F(s), f\rangle\rangle\in L^{1}(\mathbb{R},ds)$ for
 all $f\in\mathcal{S}$, then $\int_{\mathbb{R}}F(s)ds$ is defined to
 be the unique element of $(\mathcal{S})^{-1}$ such that
$$
\langle\langle \int_{\mathbb{R}}F(s)ds,
f\rangle\rangle=\int_{\mathbb{R}}\langle\langle F(s),
f\rangle\rangle ds.
$$
\textbf{Definition 3.3} Suppose that $F:\mathbb{R}\longrightarrow
 (\mathcal{S})^{-1}$
 such that $F(s)\diamond
\dot{X}^{\beta}(s) $ is $ds-$ integrable in $ (\mathcal{S})^{-1}$.
Then, we define the Skorohod integral of $F$ with respect to $
X^{\beta}$
  by
$$
\delta^{\beta}(F):=\int_{\mathbb{R}}F(s)\delta
X^{\beta}(s):=\int_{\mathbb{R}}F(s)\diamond
\dot{X}^{\beta}_{s}ds.\eqno(3.10)
$$
In particular, if $A\subset\mathbb{R}$ is a Borel set, then
$$
\int_{A}F(s)\delta
X^{\beta}(s):=\int_{\mathbb{R}}1_{A}(s)F(s)\diamond
\dot{X}^{\beta}_{s}ds.\eqno(3.11)
$$
By Definition 3.3, we get\\
 \textbf{Proposition  3.4} \ Let $F: \mathbb{R}\longrightarrow (\mathcal{S})^{-1}$ be Skorohod integrable
with respect to $ X^{\beta}$, $Y\in(\mathcal{S})^{-1}$, then
$$
Y\diamond \delta^{\beta}(F)=\delta^{\beta}(Y\diamond F),
$$
the equation holds whenever one side exists.\\
 \textbf{Proposition  3.5} \ Let $F: \mathbb{R}\longrightarrow (\mathcal{S})^{-1}$ be Skorohod integrable
with respect to $ X^{\beta}$ and $K^{\beta}F\in \mathbb{L}$, then
$$
\delta^{\beta}(F)=\delta(K^{\beta}F).
$$
\textbf{Proof:} The S-transform of $\delta^{\beta}(F)$ is given by$$
S(\int_{\mathbb{R}}F(t)\delta
X^{\beta}_{t})(\eta)=\int_{\mathbb{R}}\int_{\mathbb{R}_{0}}S(F(t))(\eta)yI_{+}^{\beta}\eta(\cdot,y)(t)\nu(dy)dt,
\eta\in U_{p}.
$$
On the other hand, by Proposition 2.8, the S-transform of
$\delta(K^{\beta}F)$ is given by
$$\aligned
S(\delta(K^{\beta}F))(\eta)&=\int_{\mathbb{R}}\int_{\mathbb{R}_{0}}S(I_{-}^{\beta}F(t)y)(\eta)\eta(t,y)\nu(dy)dt\\
&=\int_{\mathbb{R}}\int_{\mathbb{R}_{0}}S(I_{-}^{\beta}F(t)\eta(t,y))(\eta)y\nu(dy)dt\\
&=\int_{\mathbb{R}_{0}}[S(\int_{\mathbb{R}}I_{-}^{\beta}F(t)\eta(t,y)dt)(\eta)]y\nu(dy)\\
&=\int_{\mathbb{R}_{0}}[S(\int_{\mathbb{R}}F(t)I_{+}^{\beta}\eta(\cdot,y)(t)dt)(\eta)]y\nu(dy)\\
&=\int_{\mathbb{R}}\int_{\mathbb{R}_{0}}S(F(t))(\eta)yI_{+}^{\beta}\eta(\cdot,y)(t)\nu(dy)dt,\eta\in
U_{p}.\endaligned
$$
Then,
$$
S(\int_{\mathbb{R}}F(t)\delta
X^{\beta}_{t})(\eta)=S(\delta(K^{\beta}F))(\eta),
$$ Hence, we get the desired result.\hfill{$\Box$}

Especially, for $F\in L^{p}(\Omega\times\mathbb{R})$ with
$1<p<\frac{1}{\beta}$ and $I_{-}^{\beta}F\in
L^{2}(\Omega\times\mathbb{R})$, we have
$$
\int_{\mathbb{R}}F(t)\delta
X^{\beta}_{t}=\int_{\mathbb{R}}I_{-}^{\beta}F(t)dX_{t},\eqno(3.12)
$$
which holds whenever one side exists.
 \par
Let  $0<\alpha<\beta<\frac{1}{2}$, $X^{\alpha}$ is a fractional
 L\'{e}vy processes  defined by the same underlying zero mean square integrable L\'{e}vy process $X$
 of $X^{\beta}$
, that is
 $$
X_{t}^{\alpha}=\int_{\mathbb{R}}I_{-}^{\alpha}\chi_{[0,t]}(s)dX_{s}=\frac{1}{\Gamma(\alpha+1)}
\int_{-\infty}^{t}((t-s)^{\alpha}-(-s)_{+}^{\alpha})dX_{s}.
$$
Then, by the semigroup property of Riemann-Liouvile fractional
integration operator, i.e.
$I_{-}^{\alpha}I_{-}^{\beta}=I_{-}^{\alpha+\beta}$,
$0<\alpha,\beta,\alpha+\beta<1$, we get the following relation
between two different fractional L\'{e}vy processes:
$$
X^{\beta}_{t}=\int_{\mathbb{R}}I_{-}^{\alpha}I_{-}^{\beta-\alpha}1_{[0,t]}(s)dX_{s}=\frac{1}{\Gamma(\beta-\alpha)}
\int_{-\infty}^{\infty}((t-s)_{+}^{\beta-\alpha}-(-s)_{+}^{\beta-\alpha})dX^{\alpha}_{s}\eqno(3.13),
$$
which holds in $L^{2}$-sense. \\
More generally, by Definition 3.3 and the  semigroup property of
Riemann-Liouvile fractional integration operator, we get the
following integral transformation formula between two different
fractional L\'{e}vy
processes:\\
 \textbf{Proposition 3.6}~ Let
$0<\alpha<\beta<\frac{1}{2}$, $X^{\alpha}$, $X^{\beta}$  are
fractional L\'{e}vy processes given above. $F:B\times [0, T]\times
\Omega\longrightarrow (\mathcal{S})^{-1}$, $B\subset\mathbb{R}$ is a
Borel set, $F$ is  Skorohod integrable with respect to $ X^{\beta}$,
and   $K^{\beta}F\in \mathbb{L}$. Then
$K^{\alpha}I_{-}^{\beta-\alpha}F\in \mathbb{L}$ and
$$
\delta^{\beta}(F)=\delta^{\alpha}(I_{-}^{\beta-\alpha}F).\eqno(3.14)
$$
\textbf{Proof}: In fact, by the semigroup property of
Riemann-Liouvile fractional integration operator and Proposition
3.4, we have the following relations
 $$
\delta^{\beta}(F)=\delta(K^{\beta}F)=\delta(I_{-}^{\alpha}K^{\beta-\alpha}F)=\delta^{\alpha}(I_{-}^{\beta-\alpha}F).
$$
\textbf{Remark 3.7}~ Proposition 3.6 provides a powerful tool in
 solving stochastic differential equations driven by
 fractional L\'{e}vy  noises, by which one can
 transform stochastic differential equations driven by
 fractional L\'{e}vy  noises with different parameter $\beta$ into
 the same kind of stochastic differential equations with the same
 parameter.
\\[4mm]
\noindent {\bbb 4\quad The stochastic Volterra equation driven by fractional L\'{e}vy process}\\[1mm]
\noindent In this section, we consider the following Skorohod
stochastic integral equation of  Volterra-type  driven by fractional
L\'{e}vy process
$$
U(t)=a(t)+ \int_{0}^{t}b(t,s)\diamond U(s)ds+
\int_{0}^{t}\sigma(t,s)\diamond U(s)\delta X^{\beta}_{s},0\leq t\leq T. \eqno(4.1)
$$
Since (4.1) can be written as
$$
U(t)=a(t)+ \int_{0}^{t}b(t,s)\diamond U(s)ds+
\int_{0}^{t}\sigma(t,s)\diamond U(s)\diamond
\dot{X}^{\beta}_{s}ds,0\leq t\leq T.
$$
Thus, (4.1) can be regarded as a special case of the following
linear stochastic Volterra equation:
$$
U(t)=J(t)+ \int_{0}^{t}K(t,s)\diamond U(s)ds,0\leq t\leq T
\eqno(4.2)
$$
where $T> 0$ is a given number and $J: [0,T]\rightarrow
(\mathcal{S})^{-1}$, $K: [0,T]\times[0,T]\rightarrow
(\mathcal{S})^{-1}$ are given stochastic distribution processes.
Then, we first consider  the solution of (4.2) in
$(\mathcal{S})^{-1}$.
\\
\textbf{Lemma 4.1}  \ Let $J: [0,T]\rightarrow (\mathcal{S})^{-1}$,
$K: [0,T]\times[0,T]\rightarrow (\mathcal{S})^{-1}$ are continuous
stochastic distribution processes. If there exists a $q$ satisfying
$p_{0}<q<\infty$ and $M<\infty$ such that
$$
\|K(t,s)\|_{-1,-q}<M, 0\leq s\leq t\leq T,
$$
then there exists a unique continuous stochastic distribution
process solves (4.2) which is given by
$$
U(t)=J(t)+ \int_{0}^{t}H(t,s)\diamond U(s)ds, \eqno(4.3)
$$
where
$$
\aligned
& H(t,s)=\sum_{n=1}^{\infty}K_{n}(t,s),\\
& K_{n+1}(t,s)=\int_{s}^{t}K_{n}(t,u)\diamond K(u,s)du,n\geq1,\\
&K_{1}(t,s)=K(t,s).\endaligned$$
 \textbf{Proof}: The iteration method of Theorem 3.4.2 of \cite{Holden} is
 also valid in our context. Thus, we omit it here.\\
\textbf{Theorem 4.2}   \ Assume that $a : [0,T]\rightarrow
(\mathcal{S})^{-1}$ is continuous, $b: [0,T]\times[0,T]\rightarrow
(\mathcal{S})^{-1}$ and $\sigma : [0,T]\times[0,T]\rightarrow
(\mathcal{S})^{-1}$ are bounded continuous function. Then the
equation (4.1) has a unique solution in $(\mathcal{S})^{-1}$.\\
 \textbf{Proof}: Let $J(t)=a(t)$, $K(t,s)=b(t,s)+\sigma(t,s)\diamond\dot{X}^{\beta}_{s}
  $, $J(t)$ is clearly continuous. Since $b$ and
 $\sigma$ is bounded and continuous, it   suffices to investigate
 $\dot{X}^{\beta}_{t}$. By Theorem 3.1, $\dot{X}^{\beta}_{t}\in (\mathcal{S})^{-1}_{-p}$
all $p>\max\{1,q_{0}\}$ and
$$
\dot{X}^{\beta}_{t}=\langle C_{1}, \lambda_{t}\rangle,
$$
where
$$\lambda_{t}(u,y)=yI_{-}^{\beta}\delta_{t}(u)=y(t-u)_{+}^{\beta-1}/\Gamma(\beta).$$
From \cite{Thangavelu},   Hermite functions $\{\xi_{n}\}$ is the
orthogonal basis in $L^{2}(\mathbb{R})$ and
\begin{equation} |\xi_{n}(x)|\leq\begin{cases}Cn^{-\frac{1}{12}},&\  |x|\leq 2\sqrt{n},\\
Ce^{-\gamma x^{2}},&\ |x|\geq2\sqrt{n} \ ,\nonumber
\end{cases}
\end{equation}
where $C$ and $\gamma$ are certain positive constants independent of
$n$. Let $t>s$, then we have
$$ \aligned
&\langle(t-u)_{+}^{\beta-1}-(s-u)_{+}^{\beta-1},
\xi_{n}\rangle_{L^{2}(\mathbb{R})}\\
&=\int_{\mathbb{R}}[(t-u)_{+}^{\beta-1}-(s-u)_{+}^{\beta-1}]\xi_{n}(u)du\\
&=\int_{\mathbb{R}}[(t-s-u)_{+}^{\beta-1}-(-u)_{+}^{\beta-1}]\xi_{n}(s+u)du\\
&=(t-s)^{\beta}\int_{\mathbb{R}}[(1-u)_{+}^{\beta-1}-(-u)_{+}^{\beta-1}]\xi_{n}(s+(t-s)u)du\\
&\leq
(t-s)^{\beta}\{\int_{|s+(t-s)u|\leq2\sqrt{n}}Cn^{-\frac{1}{12}}[(1-u)_{+}^{\beta-1}-(-u)_{+}^{\beta-1}]du\\
& \ \ \ + \int_{|s+(t-s)u|>2\sqrt{n}}Ce^{-\gamma
u^{2}}[(1-u)_{+}^{\beta-1}-(-u)_{+}^{\beta-1}]du\}\\
&\leq
C(t-s)^{\beta}n^{\frac{1}{12}}\int_{\mathbb{R}}[(1-u)_{+}^{\beta-1}-(-u)_{+}^{\beta-1}]du\\
&\leq \widetilde{C}n^{\frac{1}{12}}(t-s)^{\beta}.
\endaligned
$$
Hence,
 $$ \aligned
&\|\dot{X}^{\beta}_{t}-\dot{X}^{\beta}_{s}\|_{-1,-p}^{2}=\|\lambda_{t}(u,y)-\lambda_{s}(u,y)\|_{-p,\pi}^{2}\\
&=(\int_{\mathbb{R}}|y|^{2}d\nu(y))\sum_{n=1}^{\infty}(n+1)^{-2p}\frac{1}{\Gamma(\beta)}\langle(t-u)_{+}^{\beta-1}-(s-u)_{+}^{\beta-1},
\xi_{n}\rangle_{L^{2}(\mathbb{R})}^{2}\\
&\leq
\widetilde{\widetilde{C}}|s-t|^{2\beta}\sum_{n=1}^{\infty}(n+1)^{-2p+\frac{1}{6}}\\
&\leq C'|s-t|^{2\beta}<\infty,  for\  p>\max\{1,q_{0}\},
\endaligned $$
where $\widetilde{\widetilde{C}}, \widetilde{C},  C' $ are positive
constants.  Thus, $\dot{X}^{\beta}_{t}$ is continuous in
$(\mathcal{S})^{-1}$. Hence, by Lemma 4.1, the equation (4.1) has a
unique solution in $(\mathcal{S})^{-1}$.\hfill{$\Box$}
\\[4mm]
\noindent {\bbb 5\quad A general existence and uniqueness theorem}\\[1mm]
\noindent In this section, we consider a   existence and uniqueness
results for the stochastic equation  satisfying linear growth  and
Lipschtz condition. First we consider a general stochastic
differential equation with Lipschtz and linear growth conditions in
the generalized distribution space, and we obtain a unique
continuous global solution in $(\mathcal{S})_{-p}^{-1}$ with
$p>q_{0}$ a natural number. Then, we consider a  stochastic
differential equation driven by fractional L\'{e}vy noise with
Lipschtz condition and linear conditions, because of the boundedness
of fractional L\'{e}vy noise in $(\mathcal{S})_{-p}^{-1}$ for
$p>\max\{1,q_{0}\}$, we obtain
  a global solution.
\\
\textbf{Lemma 5.1} Let $p>q_{0}$ be a natural number and suppose
that $F:[0, +\infty)\times (\mathcal{S})_{-p}^{-1}\rightarrow
(\mathcal{S})_{-p}^{-1}$ satisfies the following two conditions:
$$
\|F(t, Y)\|_{-1,-p}\leq C(1+\|Y\|_{-1,-p})\eqno(5.1)
$$
$$
\|F(t, Y)-F(t, Z)\|_{-1,-p}\leq C\|Y-Z\|_{-1,-p}\eqno(5.2)
$$
for all $t\in [0, +\infty)$, $Y, Z\in(\mathcal{S})_{-p}^{-1}$,  $C$
is a constant independent of $t$,  $Y$, $Z$. Then the differential
equation
$$
\frac{dU(t)}{dt}=F(t,U(t)),
U(0)=U_{0}\in(\mathcal{S})_{-p}^{-1}\eqno(5.3)
$$
has  a unique t-continuous global solution $U:[0,
+\infty)\rightarrow
(\mathcal{S})_{-p}^{-1}$.\\
 \textbf{Proof:}We can use the classical iteration methods of the linear differential equation to verify the lemma, we omit it here.
 \hfill{$\Box$}
\par  Based on the above results, we  consider the following
stochastic equation driven by fractional L\'{e}vy noise:
$$
U(t)= U_{0}+\int_{0}^{t}b(U(s))ds+ \int_{0}^{t}\sigma(U(s))
X^{\beta}_{s},t\geq 0. \eqno(5.4)
$$
which can be written as
$$
U(t)= U_{0}+\int_{0}^{t}b(U(s))ds+ \int_{0}^{t}\sigma(U(s))\diamond
\dot{X}^{\beta}_{s}ds,t\geq 0. \eqno(5.5)
$$
or
$$
\frac{dU(t)}{dt}=b(U(t))+\sigma(U(t))\diamond \dot{X}^{\beta}_{t},
U(0)=U_{0}\in(\mathcal{S})_{-p}^{-1},t\geq 0.\eqno(5.6)
$$
\textbf{Theorem 5.2}~ Let $p>\max\{1,q_{0}\}$ be a natural number,
suppose that $b:(\mathcal{S})_{-p}^{-1}\rightarrow
(\mathcal{S})_{-p}^{-1}$ and
$\sigma:(\mathcal{S})_{-p}^{-1}\rightarrow (\mathcal{S})_{-p}^{-1}$
satisfies the following  conditions:
$$\aligned
& \|b(Y)\|_{-1,-p}\leq C(1+\|Y\|_{-1,-p}),\\
&\|b(Y)-b(Z)\|_{-1,-p}\leq C\|Y-Z\|_{-1,-p},\\
& \|\sigma(Y)\|_{-1,-p}\leq C(1+\|Y\|_{-1,-p}),\\
&\|\sigma(Y)-\sigma(Z)\|_{-1,-p}\leq C\|Y-Z\|_{-1,-p},
\endaligned$$
for all   $Y, Z\in(\mathcal{S})_{-p}^{-1}$, with $C$ independent of
   $Y$, $Z$. Then the differential equation (5.5) has a unique
continuous solution $U:[0, \infty)\rightarrow
(\mathcal{S})_{-p}^{-1}$.\\
\textbf{Proof:} Let $F(t,Y)=b(Y)+\sigma(Y)\diamond
\dot{X}^{\beta}_{t}$, by (3.8) the  deduction in Theorem 4.2, we can
get for $p>\max\{1,q_{0}\}$, $\forall t\geq
0,\|\dot{X}^{\beta}(t)\|_{-1,p}\leq M $, then
$$\aligned
\|F(t, Y)\|_{-1,-p}&=\|b(Y)+\sigma(Y)\diamond
\dot{X}^{\beta}_{t}\|_{-1,-p}\\
&\leq \|b(Y)\|_{-1,-p} + \|\sigma(Y)\diamond
\dot{X}^{\beta}_{t}\|_{-1,-p} \\
&=\|b(Y)\|_{-1,-p} + \|\sigma(Y)\|_{-1,-p}
\|\dot{X}^{\beta}(t)\|_{-1,-p}\\
&\leq C(1+M)(1+\|Y\|_{-1,-p})\endaligned
$$
Similarly,
$$
\|F(t, Y)-F(t, Z)\|_{-1,-p}\leq C(1+M)\|Y-Z\|_{-1,-p}.
$$
That is, $F$ satisfies   Lipschtz condition and linear growth
condition. Hence, by  Lemma 5.1, we deduce
that  the stochastic equation (5.5) has a unique global solution
$U:[0, \infty)\rightarrow (\mathcal{S})_{-p}^{-1}$.\hfill{$\Box$}
\noindent \\[4mm]
 \noindent{\bbb{References}}


\begin{enumerate}
{\footnotesize \bibitem{Lokka}Lokka A,   Proske F. Infinite dimensional analysis of pure jump L\'{e}vy processes on the Poisson space.
 Mathematica Scandinavica. 2006, 98: 237-261\\[-6mm]

\bibitem{Kol1}Kolmogorov A N. Wienersche Spiralen und einige andere interessante Kurven
          in Hilbertschen Raum.   C  R  (Doklady) Acad  Sci USSR (NS). 1940, 26: 115-118 \\[-6mm]

\bibitem{Mandelbrot}Mandelbrot B, Van Vess J. Fractional Brownian motion,
fractional noises and application. SIAM Rev. 1968, 10: 427-437 \\[-6mm]

\bibitem{Marquardt} Marquardt T.  Fractional L\'{e}vy processes with an
application to long memory moving average processes. Bernoulli. 2006, 12: 1099-1126 \\[-6mm]

\bibitem{HuangLC} Huang Z, Li C.  On fractional stable processes and sheets: white noise approach. J
 Math  Anal  Appl.  2007, 325: 624-635 \\[-6mm]

\bibitem{HuangLP} Huang Z, Li P. Generalized  fractional L\'{e}vy
processes: a white noise approach.  Stoch  Dyn. 2006, 6: 473-485
\\[-6mm]
\bibitem{HuangLU} Huang Z, L\"{u} X, Wan J.  Fractional L\'{e}vy
processes and noises
on Gel$\prime$fand triple.  Stoch Dyn.   2010, 10: 37-51 \\[-6mm]
\bibitem {LUHuang}  L\"{u} X, Huang Z, Wan J.  Fractional L\'{e}vy
Processes on
  and Stochastic Integration.  Front Math
 China. 2008, 3: 287-303  \\[-6mm]
\bibitem{Bender} Bender C, Marquardt T.  Stochastic calculus for convoluted  L\'{e}vy processes. Bernoulli. 2008, 14(2): 499--518 \\[-6mm]

\bibitem{Sato} Sato K. L\'{e}vy Processes and Infinitely Divisible
Distributions. Cambridge: Cambridge University Press, 1999\\[-6mm]

\bibitem{Samko} Samko S G ,  Kilbas A A,  Marichev O I. Fractional Integrals
and Derivatives: Theory and Applications. Gordon and Breach, 1987\\[-6mm]

\bibitem {Elliott} Elliott R C, Van der Hoek J. A general fractional white noise
theory and applications  to finance.  Mathematical Finance. 2003, 13
:
301-330\\[-6mm]

\bibitem{Holden} Holden H, Oksendal B, Uboe J, Zhang T. Stochastic Partial
Differential Equations: A modeling, white noise functional approach.
Birkhauser, 1996\\[-6mm]

\bibitem{Thangavelu} Thangavelu S, Lectures of Hermite and Laguerr Expansions. Princeton University Press, 1993}\end{enumerate}
\end{document}